\documentclass[final]{siamltex}
\pdfoutput=1
\usepackage{amssymb,amsfonts,amsmath}
\usepackage[mathscr]{eucal}
\usepackage{stmaryrd}
\usepackage{amsbsy}
\usepackage{bm}
\usepackage{braket}
\usepackage[usenames]{xcolor}
\definecolor{color1}{rgb}{0.2235,0.4157,0.6941}
\definecolor{color2}{rgb}{0.8549,0.4863,0.1882}
\definecolor{color3}{rgb}{0.2431,0.5882,0.3176}
\usepackage{array}
\usepackage{tabularx}
\usepackage[colorlinks,urlcolor=color3,citecolor=color3,linkcolor=color3]{hyperref}
\usepackage{algorithm}
\usepackage{algpseudocode}
\usepackage{mathtools}
\usepackage{tikz}
\usetikzlibrary{calc}

\newcommand{\Parens}[1]{\left(#1\right)}
\newcommand{\Tra}{{\sf T}} 
\newcommand{\Real}{\mathbb{R}}
\newcommand{\qtext}[1]{\quad\text{#1}\quad}
\newcommand{\V}[2][]{{\bm{#1\mathbf{\MakeLowercase{#2}}}}} 
\newcommand{\M}[2][]{{\bm{#1\mathbf{\MakeUppercase{#2}}}}} 
\newcommand{\MC}[3][]{\V[#1]{#2}_{#3}} 
\newcommand{\T}[2][]{\boldsymbol{#1\mathscr{\MakeUppercase{#2}}}} 
\newcommand{\TE}[3][]{#1{\MakeLowercase{#2}}_{#3}}

\newcommand{\Eqn}[1]{\hyperref[eq:#1]{{\rm (\ref*{eq:#1})}}} %
\newcommand{\Fig}[1]{\hyperref[fig:#1]{Figure~\ref*{fig:#1}}} %
\newcommand{\Tab}[1]{\hyperref[tab:#1]{Table~\ref*{tab:#1}}} %
\newcommand{\Alg}[1]{\hyperref[alg:#1]{Algorithm~\ref*{alg:#1}}} %

\begin{document}

\title{Symmetric Orthogonal Tensor Decomposition is Trivial\thanks{This material is based upon work supported by the U.S. Department of
Energy, Office of Science, Office of Advanced Scientific Computing
Research, Applied Mathematics program.
Sandia National Laboratories is a multi-program laboratory managed and
operated by Sandia Corporation, a wholly owned subsidiary of Lockheed
Martin Corporation, for the U.S. Department of Energy's National
Nuclear Security Administration under contract DE--AC04--94AL85000. 
}}
\author{%
  Tamara G. Kolda\footnotemark[2] 
}
\date{}
\maketitle

\renewcommand{\thefootnote}{\fnsymbol{footnote}}
\footnotetext[2]{Sandia National Laboratories, Livermore, CA.
  Email: tgkolda@sandia.gov}
\renewcommand{\thefootnote}{\arabic{footnote}}

\begin{abstract}
  We consider the problem of decomposing a real-valued symmetric
  tensor as the sum of outer products of real-valued, pairwise
  orthogonal vectors. Such decompositions do not generally exist, but
  we show that some symmetric tensor decomposition problems can be
  converted to orthogonal problems following the whitening procedure
  proposed by Anandkumar et~al.~(2012).  If an orthogonal
  decomposition of an $m$-way $n$-dimensional symmetric tensor
  exists, we propose a novel method to compute it that reduces to an
  $n \times n$ symmetric matrix eigenproblem. We provide numerical
  results demonstrating the effectiveness of the method.
\end{abstract}

\section{Introduction}
\label{sec:intro}
Let $\T{A}$ be an $m$-way $n$-dimensional real-valued symmetric tensor. 
Let $\T{A}$ represent an $m$-way, $n$-dimension symmetric tensor. 
Given a real-valued vector $\V{x}$ of length $n$, 
we let $\V{x}^m$ denote the $m$-way, $n$-dimensional 
symmetric outer product tensor such that 
$\Parens{\V{x}^m}_{i_1 i_2 \cdots i_m} = x_{i_1} x_{i_2} \cdots x_{i_m}$.
Comon et~al.~\cite{CoGoLiMo08} showed there exists a decomposition of the form
\begin{equation}
  \label{eq:scp}
  \T{A} = \sum_{k=1}^p \lambda_k \MC{X}{k}^m,
\end{equation}
where $\V{\lambda} = [\lambda_1 \cdots \lambda_p]^{\Tra} \in \Real^p$
and $\M{X} = [\MC{x}{1} \cdots \MC{x}{p}] \in \Real^{n \times p}$;
see \Fig{symcp}.
Without loss of generality, we assume each $\MC{x}{k}$ has unit norm, i.e., $\|\MC{x}{k}\|_2 = 1$.
The least value $p$ such that \Eqn{scp} holds is called the symmetric tensor rank.
Finding real-valued symmetric tensor decompositions has been the topic of several recent papers, e.g., \cite{SymCP-arXiv-1410.4536}.

\definecolor{TikzCubeColor}{rgb}{0.7,0.7,0.7}
\newcommand{\TikzCube}[3]{ %
  \def\scale{#1};
  \def\xsize{1*\scale};
  \def\ysize{1*\scale};
  \def\xdelta{0.40*\scale};
  \def\ydelta{0.25*\scale};

  \coordinate (FrontLowerLeft) at (#2);
  \coordinate (TopLowerLeft) at ($(FrontLowerLeft)+(0,\ysize)$);
  \coordinate (SideLowerLeft) at ($(FrontLowerLeft)+(\xsize,0)$);
  \coordinate (FrontCenter) at ($(FrontLowerLeft)+(\xsize/2,\ysize/2)$);

  \draw[fill=TikzCubeColor!50] (FrontLowerLeft) -- ++(\xsize,0) -- ++(0,\ysize) -- ++(-\xsize,0) -- cycle;
  \draw[fill=TikzCubeColor!50] (TopLowerLeft)  -- ++(\xdelta,\ydelta) -- ++(\xsize,0) -- ++(-\xdelta,-\ydelta) -- cycle;
  \draw[fill=TikzCubeColor!50] (SideLowerLeft) -- ++(\xdelta,\ydelta) -- ++(0,\ysize) -- ++(-\xdelta,-\ydelta) -- cycle;

  \node at (FrontCenter) {#3};
}
\newcommand{\TikzStar}[3]{
  \def\scale{#1};
  \def\xsize{1*\scale};
  \def\ysize{1*\scale};
  \def\xdelta{0.40*\scale};
  \def\ydelta{0.25*\scale};
  \def\exdelta{0.15*\scale};
  \def\width{0.1*\scale};

  \coordinate (FrontUpperLeft) at ($(#2)+(0,\ysize)$);
  \coordinate (ColUpperLeft) at ($(FrontUpperLeft)+(0,-\exdelta/4)$);
  \coordinate (TubeLowerLeft) at ($(FrontUpperLeft)+(0,\exdelta/4)$);
  \coordinate (RowUpperLeft) at ($(FrontUpperLeft)+(\exdelta,0)$);

  \draw[fill=TikzCubeColor!50] (RowUpperLeft) node[left=2] {$\lambda_{#3}$} -- ++(0,-\width) -- ++(\xsize,0) -- ++(0,\width) node[right] {$\MC{x}{#3}$} -- cycle;
  \draw[fill=TikzCubeColor!50] (TubeLowerLeft) -- ++(1.25*\width,0) -- ++(\xdelta,\ydelta) -- ++(-1.25*\width,0) node[above] {$\MC{x}{#3}$} -- cycle;
  \draw[fill=TikzCubeColor!50] (ColUpperLeft) -- ++(\width,0) -- ++(0,-\ysize) -- ++(-\width,0) node[right=2] {$\MC{x}{#3}$} -- cycle;
 
}

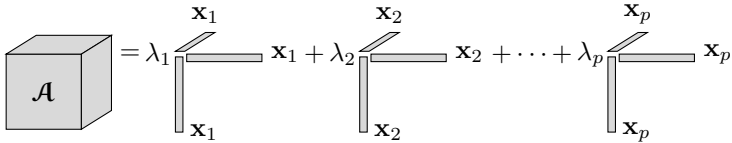
\begin{figure}[htbp]
  \centering
  \begin{tikzpicture}
    \TikzCube{1}{0,0}{$\T{A}$};
    \node at (1.65,1) {$=$};
    \TikzStar{1}{2.25,0}{1};
    \node at (4.1,1) {$+$};
    \TikzStar{1}{4.7,0}{2};
    \node[right] at (6.33,1) {$+\cdots+$};
    \TikzStar{1}{8,0}{p};   
  \end{tikzpicture}
  \caption{Symmetric tensor factorization for $m=3$.}
  \label{fig:symcp}
\end{figure}

If we can discover an $\M{X}$ with orthogonal columns, i.e.,
$\M{X}^{\Tra}\M{X} = \M{I}_p$, then we say that $\T{A}$ has an
orthogonal symmetric tensor decomposition.  Generally, orthogonal
decompositions \emph{do not} exist; Robeva~\cite{Ro14} classifies the
tensors that have such decompositions.
Nevertheless, we consider the problem of how to compute orthogonal
symmetric tensor decompositions, as has been recently considered by
Anandkumar et~al.~\cite{AnGeHsKa12}.
They show that the pairs $(\lambda_k, \MC{x}{k})$ are Z-eigenpairs of
$\T{A}$ for $k=1,\dots,p$, and propose solving the problem via an
iterative power method.
We show that this problem can instead be solved via a symmetric matrix
eigenproblem on an $n \times n$ matrix, including determining the
rank.  This is can be interpreted as a special case of the
simultaneous matrix diagonalization approach proposed by
De~Lathauwer~\cite{De06}.

Although orthogonal symmetric tensor decompositions do not generally
exist, Anandkumar et~al.~\cite{AnGeHsKa12} showed that certain
symmetric tensor decompositions problems can be transformed via \emph{whitening} to orthogonal problems. We generalize their results to show that such a transformation is possible whenever $\M{X}$ has full column rank and there exists a linear combination of two-dimensional slices of $\T{A}$ that is positive definite. The second condition is always satisfied, for example, if $\T{A}$ is positive definite.

\section{Background}

A tensor is a multidimensional array. The number of ways or modes is
called the \emph{order} of a tensor. For example, a matrix is a tensor
of order two. Tensors of order three or greater are called
\emph{higher-order} tensors.

We use the notation $\M{I}_n$ to denote the $n \times n$ identity matrix.

We use the term \emph{generic} to mean with probability one.
For instance, it is well known that a random $n \times n$ matrix generically has rank $n$. 

\subsection{Symmetry}
A tensor is symmetric if its entries do not change under permutation of the indices.
Formally, 
we let $\pi(m)$ denote the set of permutations of length $m$. 
We say a real-valued $m$-way $n$-dimensional tensor 
$\T{A}$ is \emph{symmetric} \cite{CoGoLiMo08} if
\begin{displaymath}
  \TE{a}{i_{p(1)} \cdots i_{p(m)}} = \TE{a}{i_1 \cdots i_m}
  \qtext{for all} i_1, \dots, i_m \in \set{1,\dots,n}
  \text{ and } p \in \pi(m).
\end{displaymath}

\subsection{Tensor-vector Products}
The tensor-vector product $\T{A}\V{x}^{(m-1)}$ produces a vector in $\Real^n$ such that 
\begin{displaymath}
  \Parens{ \T{A} \V{x}^{m-1} }_{i_1} = \sum_{i_2,\dots,i_m=1}^n a_{i_1\cdots i_m} x_{i_2} \cdots x_{i_m}
  \qtext{for}
  i_1 \in \set{1,\dots,n}.
\end{displaymath}
The tensor-vector product $\T{A}\V{x}^m$ produces a scalar such that 
\begin{displaymath}
  \T{A}\V{x}^m = \V{x}^{\Tra} \Parens{ \T{A} \V{x}^{m-1} } = \sum_{i_1,\dots,i_m=1}^n a_{i_1\dots i_m} x_{i_1} \cdots x_{i_m}.
\end{displaymath}
A symmetric tensor $\T{A}$ is positive definite if $\T{A}\V{x}^m > 0$ for all $\V{x} \neq 0$.

\subsection{Tensor-matrix Products}
Let $\M{V}$ be an $p \times n$ matrix. Then the tensor-matrix product $\T{A}(\M{V},\dots,\M{V})$ indicates multiplication of the tensor $\T{A}$ in each mode by the matrix $\M{V}$. The result is a symmetric tensor that is the same order as $\T{A}$ but now $p$-dimensional such that
\begin{displaymath}
  \Parens{\T{A}(\M{V},\dots,\M{V})}_{i_1 \dots i_m} = \sum_{j_1,\dots,j_m=1}^n a_{j_1 \dots j_m} v_{i_1 j_1} \cdots v_{i_m j_m}
  \qtext{for} i_1,\dots,i_m \in \set{1,\dots,p}.
\end{displaymath}

\subsection{Tensor Rank}      
Recall that the rank of a tensor $\T{A}$, denoted $\text{rank}(\T{A})$,
is the smallest number of rank-one tensors that sums to the original
tensor \cite{KoBa09}.
The symmetric tensor rank, denoted $\text{symrank}(\T{A})$, is the
smallest number of symmetric rank-one tensors that sums to the
original tensor \cite{CoGoLiMo08}.
In the case of the tensor $\T{A}$ in \Eqn{scp} with $\M{X}$ having orthogonal columns, it is easy to show that the symmetric tensor rank of $\T{A}$ is equal to the tensor rank of $\T{A}$ which is equal to the matrix rank of $\M{X}$ which is equal to $p$, i.e., 
\begin{displaymath}
  \text{symrank}(\T{A}) = \text{rank}(\T{A}) = \text{rank}(\M{X}) = p,
\end{displaymath}
via unfolding arguments.

\section{Orthogonal Symmetric Decomposition}
Given $\T{A}$, our goal is to find $\V{\lambda}$ and $\M{X}$ that
satisfies \Eqn{scp}, under the assumption that $\M{X}$ is known to
have orthogonal columns.

We define $\M{B}$ to be an arbitrary linear combination
of slices of $\T{A}$:
\begin{equation}
  \label{eq:bgen}
  \M{B} =  \sum_{i_3,\dots,i_m} \beta_{i_3 \cdots i_m} \M{A}(:,:,i_3,\dots,i_m) 
  = \sum_{k=1}^p \sigma_k \MC{x}{k}\MC{x}{k}^{\Tra}
  = \M{X} \M{\Sigma} \M{X}^{\Tra}
\end{equation}
where the $\beta$-values define the linear combination and 
\begin{displaymath} 
  \sigma_k =  \lambda_k \sum_{i_3,\dots,i_m} \beta_{i_3 \cdots i_m}
  x_{i_3 k} \cdots x_{i_m k}
  \qtext{for} k = 1,\dots,p.
\end{displaymath}
Observe that $(\sigma_k,\V{x_k})$ is an eigenpair of $\M{B}$.

\subsection{Generic case}
If $\M{X}$ is a generic matrix (with orthogonal columns) and the $\beta$-values are
arbitrary, then $\lambda_k \neq 0$ generically implies $\sigma_k \neq 0$ for
$k=1,\dots,p$. It follows that the tensor rank of $\T{A}$ is generically equal to the matrix rank of $\M{B}$, i.e.,
\begin{displaymath}
  \text{rank}(\T{A}) = \text{rank}(\M{B}).
\end{displaymath}

Additionally, the $\sigma$-values are generically distinct, so the
eigenvectors of $\M{B}$ are unambiguous and equal to the vectors in
the decomposition of $\T{A}$.  Therefore, we can compute $\M{X}$ from
$\M{B}$ via an eigenvalue decomposition and recover the
$\lambda$-values via
\begin{equation}
  \label{eq:evals}
  \lambda_k = \T{A} \MC{x}{k}^m 
  \qtext{for}
  k = 1,\dots,p.
\end{equation}

\subsection{Non-generic Case}
In the generic case, the matrix rank of $\M{B}$ is equal to the symmetric
tensor rank of $\T{A}$ and the nonzero eigenvalues of $\M{B}$ are distinct.
In order to guard against the non-generic case, randomly project the
tensor as follows. Let $\M{V} \in \Real^{n \times n}$ be a random
orthonormal matrix. Then compute,
\begin{displaymath}
  \T[\hat]{A} = \T{A}(\M{V},\dots,\M{V}) 
  = \sum_{k=1}^p \lambda_k (\M{V}\MC{x}{k})^m 
  = \sum_{k=1}^p \lambda_k \MC[\hat]{x}{k}^m .
\end{displaymath}
Apply the procedure outlined above to obtain the decomposition of
$\T[\hat]{A}$ in terms of $\V[\hat]{\lambda}$ and $\M[\hat]{X}$.
Then the decomposition of $\T{A}$ is given by 
\begin{displaymath}
\V{\lambda} = \V[\hat]{\lambda}
\qtext{and}
\M{X} = \M{V}^{\Tra}\M[\hat]{X}.
\end{displaymath}

\subsection{Algorithm}
The algorithm for computing the orthogonal symmetric decomposition of
$\T{A}$ using optional randomization is given in \Alg{ostd}.

\begin{algorithm}[htpb]
  \caption{Orthogonal Symmetric Decomposition}
  \label{alg:ostd}
  Input: Let $\T{A}$ be a symmetric $m$-way, $n$-dimensional real-valued tensor
  that is known to have an orthogonal symmetric tensor decomposition.
  \begin{algorithmic}[1]
    \If{apply randomization}
    \State $\M{V} \gets$ random $n \times n$ orthonormal matrix
    \Else
    \State $\M{V} \gets \M{I}_n$
    \EndIf
    \State $\T[\hat]{A} \gets \T{A}(\M{V},\dots,\M{V})$
    \State $\beta \gets$ arbitrary $(m-2)$-way real-valued tensor of dimension $n$
    \State $\M{B} \gets \sum_{i_3,\dots,i_m} \beta_{i_3 \cdots i_m} \M[\hat]{A}(:,:,i_3,\dots,i_m)$
    \State $\set{\sigma_k, \MC[\hat]{X}{k}}_{k=1}^p \gets$  eigenpairs of $\M{B}$ with $\sigma_k \neq 0$
    \For{$k=1,\dots,p$}
    \State $\MC{x}{k} \gets \M{V}^{\Tra}\MC[\hat]{x}{k}$ %
    \State $\lambda_k \gets \T{A}\MC{x}{k}^m$ %
    \EndFor
  \end{algorithmic}
\end{algorithm}

\section{Whitening}
Although most tensors do not have symmetric decompositions, Anandkumar
et~al.~\cite{AnGeHsKa12} show how whitening may be used in a special case
of nonorthogonal symmetric tensor decomposition. We generalize their result.

Given $\T{A}$, our goal is to find $\V{\lambda}$ and $\M{X}$ that
satisfies \Eqn{scp}, under the assumption that $\M{X}$ is known to
have full column rank but may not have orthogonal columns.

\subsection{Transformation to orthogonal problem}
Let $\M{C}$ be an arbitrary linear combination
of slices:
\begin{equation}
  \label{eq:cgen}
  \M{C} =  \sum_{i_3,\dots,i_m} \gamma_{i_3 \cdots i_m} \M{A}(:,:,i_3,\dots,i_m),
\end{equation}
where the $\gamma$-values define the linear combination.
If $\M{C}$ is positive semi-definite (p.s.d.) and has the same rank as $\T{A}$, then we can apply whitening
as follows. Let 
\begin{displaymath}
  \M{U} \M{D} \M{U}^{\Tra} = \M{C}
\end{displaymath}
be the ``skinny'' eigendecomposition of $\M{C}$ where $\M{U}$ is an orthogonal
matrix of size $n \times p$ and $\M{D}$ is a diagonal matrix of size $p
\times p$.
Define the whitening matrix as
\begin{displaymath}
  \M{W} = \M{D}^{-1/2}\M{U}^{\Tra}
  \qtext{so that}
  \M{W}\M{C}\M{W}^{\Tra} = (\M{W}\M{U}\M{D}^{1/2})(\M{W}\M{U}\M{D}^{1/2})^{\Tra} = \M{I}_p.
\end{displaymath}

We use $\M{W}$ to whiten the tensor as
\begin{displaymath}
  \T[\bar]{A} = \T{A}(\M{W}, \dots, \M{W})
  = \sum_{k=1}^p \lambda_k \M{W} \MC{x}{k}
  = \sum_{k=1}^p \lambda_k \MC[\bar]{x}{k}
\end{displaymath}
Now, $\M[\bar]{X} = \M{W}\M{X} \in \Real^{p \times p}$ is a matrix with orthogonal columns. Moreover, the size of the problem is reduced because
$\T[\hat]{A}$ is an $m$-way $p$-dimensional tensor.

We compute the orthogonal tensor decomposition of $\T[\hat]{A}$ via
the procedure outlined above to get $\V[\bar]{\lambda}$ and
$\M[\bar]{X}$.
The final decomposition is given by
\begin{displaymath}
  \V{\lambda} = \V[\bar]{\lambda} 
  \qtext{and}
  \M{X} = \M{W}^{\dagger} \M[\bar]{X}.
\end{displaymath}
Here $\M{W}^{\dagger} = \M{U}\M{D}^{1/2}$ represents the psuedoinverse of $\M{W}$.

\subsection{Failure of the Method}
If the algorithm cannot find a p.s.d.\@ matrix $\M{C}$,
then the algorithm fails. Additionally, if $\text{rank}(\M{C}) < p$,
then the algorithm has a soft failure.

\subsection{Algorithm}
The algorithm for computing the symmetric decomposition of $\T{A}$
using whitening is given in \Alg{wostd}.

\begin{algorithm}[htpb]
  \caption{Whitening for Orthogonal Symmetric Decomposition}
  \label{alg:wostd}
  Input: Let $\T{A}$ be an $m$-way, $n$-dimensional real-valued tensor.
  \begin{algorithmic}[1]
    \If{apply random orthogonal matrix}
    \State $\M{V} \gets$ random $n \times n$ orthonormal matrix
    \Else
    \State $\M{V} \gets \M{I}_n$
    \EndIf
    \State $\T[\hat]{A} \gets \T{A}(\M{V},\dots,\M{V})$
    \Repeat
    \State $\gamma \gets$ arbitrary $(m-2)$-way real-valued tensor of
    dimension $n$
    \State $\M{C} \gets \sum_{i_3,\dots,i_m} \gamma_{i_3 \cdots i_m} \M[\hat]{A}(:,:,i_3,\dots,i_m)$
    \Until{$\M{C}$ is p.s.d.\@ or exit with failure}
    \State $\M{U} \M{D} \M{U}^{\Tra} \gets$ ``skinny'' eigendecomposition of $\M{C}$
    \State $\M{W} \gets \M{D}^{-1/2} \M{U}^{\Tra}$
    \State $\T[\bar]{A} \gets \T{A}(\M{W}, \dots, \M{W})$
    \State $\beta \gets$ arbitrary $(m-2)$-way real-valued tensor of same dimension as $\T[\hat]{A}$
    \State $\M{B} \gets \sum_{i_3,\dots,i_m} \beta_{i_3 \cdots i_m} \M[\bar]{A}(:,:,i_3,\dots,i_m)$
    \State $\set{\sigma_k, \MC[\bar]{X}{k}}_{k=1}^p \gets$  eigenpairs of $\M{B}$ with $\sigma_k \neq 0$
    \For{$k=1,\dots,p$}
    \State $\MC{x}{k} \gets \M{V}^{\Tra}\M{U}\M{D}^{1/2}\MC[\bar]{x}{k}$ %
    \State $\lambda_k \gets \T[\bar]{A}\MC[\bar]{x}{k}^m$ %
    \EndFor
  \end{algorithmic}
\end{algorithm}

\section{Numerical Results}
We consider the results of applying the algorithms to numerical
examples.
All experiments are done in MATLAB, Version R2014b.
Numerically, we say 
\begin{itemize}
\item an eigenvalue $\sigma$ of $\M{B}$ is nonzero if
$|\sigma| > 10^{-10}$,
\item  $\M{C}$ is p.s.d.\@ if its
smallest eigenvalue satisfies $d > -10^{-10}$,  and
\item the
skinny decomposition of $\M{C}$ uses only eigenvalues such
that $d > 10^{-10}$.
\end{itemize}
We choose both $\beta$ and $\gamma$ values from $U[0,1]$, i.e.,
uniform random on the interval $[0,1]$ and then normalize so the
values sum to one.
Random orthogonal matrices are generated via the MATLAB code
\texttt{RANDORTHMAT} by Olef Shilon.

We generate artificial data as follows.
For a given $\V{\lambda}^* \in \Real^p$ and
$\M{X}^* \in \Real^{n \times p}$, the noise-free data tensor is given by
\begin{equation}
  \label{eq:Astar}
  \T{A}^* = \sum_{k=1}^p \lambda_k^* (\MC{x}{k}^* )^m.
\end{equation}
The data tensor $\T{A}$ may also be contaminated
by noise as controlled by the parameter $\eta \geq 0$, i.e.,
\begin{equation}
  \label{eq:noise}
  \T{A} = \T{A}^*  + \eta \frac{\|\T{A}^*\|}{\|\T{N}\|} \T{N}
  \qtext{where}
  n_{i_1,\dots,i_m} \sim \mathcal{N}(0,1).
\end{equation}
Here $\T{N}$ is a noise tensor such that each element is drawn from a
normal distribution, i.e., $n_{i_1,\dots,i_m} \sim \mathcal{N}(0,1)$.
The parameters $m$, $n$, $p$ control the size of the problem. 

In our randomized experiments, we consider three sizes:
\begin{itemize}
\item $m=3,n=4,p=2$;
\item  $m=4,n=25,p=3$; and
\item $m=6,n=6,p=4$.
\end{itemize}
For each size, we also consider two noise levels: $\eta \in
\set{0,0.01}$, i.e., no noise and a small amount of noise.

The output of each run is a rank $p$, a weight vector $\V{\lambda}$, and a matrix $\M{X}$.
The relative error measures the proportion of the
observed data that is explained by the model, i.e., 
\begin{displaymath}
  \text{relative error } = 
  {\left \| \T{A} - \displaystyle\sum_{k=1}^p \lambda_k \MC{x}{k}^m  \right\|} /
  { \| \T{A} \|}.
\end{displaymath}
In the case of no noise, the ideal relative error is zero; otherwise,
we hope for something near the noise level, i.e., $\eta$.

To compare the recovered solution $\V{\lambda}$ and $\M{X}$ with the
true solution $\V{\lambda}^*$ and $\M{X}^*$, we compute the solution
score as follows.
Without loss of generality, we assume
both $\M{X}$ and $\M{X}^*$ have normalized columns. (If $\| \MC{x}{k} \|_2 \neq
1$, then we rescale $\lambda_k = \lambda_k \sqrt[m]{ \| \MC{x}{k} \|}$ and $\MC{x}{k} =
\MC{x}{k} / \| \MC{x}{k} \|$.)  There is a permutation ambiguity, but we permute
the computed solution so as to maximize the following score:
\begin{displaymath}
  \text{solution score } = \frac{1}{p} \sum_{k=1}^p 
  \left(
    1 - \frac{|\lambda_k - \lambda_k^*|}{ \max\{|\lambda_k|,|\lambda_k^*|\} }
  \right)
  \left| 
    \MC{x}{k}^{\Tra} \MC{x}{k}^*
  \right|.
\end{displaymath}
A solution score of 1 indicates a perfect match.
If $\M{X}$ has more columns than $\M{X}^*$, we choose the $p$ columns
that maximize the score.

\subsection{Orthogonal Example Showing Impact of $\beta$-values}

We discuss why we recommend a linear combination of slices instead of a
single slice.
Consider the following example. Let $m=3$, $n=3$, and $p=3$. Further,
supposed $\M{X}^* = \M{I}_n$ and $\V{\lambda}^*$ is an arbitrary vector.
Assume no noise so that $\T{A}=\T{A}^*$.
Each slide of the tensor $\T{A}$ is a rank-1 matrix. For instance,
\begin{displaymath}
  \M{A}(:,:,1) = 
  \begin{bmatrix*}
    \lambda_1 & 0 & 0 \\
    0 & 0 & 0 \\
    0 & 0 & 0
  \end{bmatrix*}.
\end{displaymath}
If we only select the first slice, which corresponds to $\beta = [1\;
0 \; 0]$, will not yield a matrix $\M{B}$ that has rank equal to $p$.
But, using a random linear combination remedies this
problem. Alternatively, using the randomization to make the problem
generic will also correct the problem.

\subsection{Random Orthogonal Examples}
In this case, we generate tensors such that $\M{X}^*$ is a random
orthogonal matrix and $\M{\lambda}^*$ is the all ones vector.
(Note that a matrix that had repeated eigenvalues would not have a
unique factorization, but tensors with repeated eigenvalues are not
impacted in the same way.)
We apply \Alg{ostd} with no randomization ($\M{V}=\M{I}_n$).
The results are shown in \Tab{orthog}.
We generated 100 random instances for each size, and then we ran the
code 10 timers per instance (each run uses a different random choice
for $\beta$).
In the noise-free case ($\eta=0$), the method works perfectly: the
rank is perfectly predicted and the exact solution is found.
In the noisy case ($\eta=0.1$), the method is less reliable. The rank
is never predicted correctly; instead, the $\M{B}$ matrix is nearly
always full rank.
Nevertheless, the relative error is usually less than $10\eta$ and the
solution score is nearly always $\geq 0.99$. So, we can presumably
threshold the small $\lambda$-values in the noisy cases to recover a
good solution.

\begin{table}[htbp]
  \centering
    \begin{tabular}{|c@{\,\,\,}c@{\,\,\,}c| %
        >{\color{color1}}r@{~~}  %
        >{\color{color2}}r@{~~}
        >{\color{color3}}r|
        >{\color{color1}}r@{~~}  %
        >{\color{color2}}r@{~~}
        >{\color{color3}}r|
        >{\color{color1}}r@{~~}  %
        >{\color{color2}}r@{~~}
        >{\color{color3}}r|}
    \hline
    \multicolumn{3}{|c}{Size} & 
    \multicolumn{9}{|c|}{No noise $\eta=0$} \\
    $m$ & $n$ & $p$ &
    \multicolumn{3}{c|}{Rank $=p$} &
    \multicolumn{3}{c|}{Rel.~Error $\leq$ 1e-10} & 
    \multicolumn{3}{c|}{Soln.~Score $\geq 0.99$} \\ 
    \hline
3 &  4 & 2 & 1000 & 100 & 2.000 & 1000 & 100 & 0.0000 & 1000 & 100 & 1.0000 \\
4 & 25 & 3 & 1000 & 100 & 3.000 & 1000 & 100 & 0.0000 & 1000 & 100 & 1.0000 \\
6 &  6 & 4 & 1000 & 100 & 4.000 & 1000 & 100 & 0.0000 & 1000 & 100 & 1.0000 \\
    \hline
    \multicolumn{3}{|c}{Size} & 
    \multicolumn{9}{|c|}{Noise $\eta=0.01$} \\
    $m$ & $n$ & $p$ &
    \multicolumn{3}{c|}{Rank $=p$} &
    \multicolumn{3}{c|}{Rel.~Error $\leq$ 0.1} & 
    \multicolumn{3}{c|}{Soln.~Score $\geq 0.99$} \\ 
    \hline
3 &  4 & 2 & 0 & 0 &  4.000 &  921 & 100 & 0.0428 &  925 & 100 & 0.9903 \\
4 & 25 & 3 & 0 & 0 & 24.999 &  818 & 100 & 0.0716 &  864 & 100 & 0.9856 \\
6 &  6 & 4 & 0 & 0 &  6.000 &  393 &  92 & 0.1916 &  498 &  99 & 0.9530 \\
\hline
    \end{tabular}
  \caption{Random orthogonal examples. We create 100 random instances
    and run the method 10 times per instance for a total of 1000
    runs. For each metric, we report \textcolor{color1}{the total number of
      runs where the metric meets the desired criteria},
    \textcolor{color2}{the number of instances where at least one run
      meets the desired criteria}, and \textcolor{color3}{the mean
      value of the metric}.}
  \label{tab:orthog}
\end{table}

\subsection{Non-orthogonal Example}

In Example 5.5(i) of \cite{Ni14}, Nie considers an method for
determining the rank of a tensor. The example tensor is of order $m=4$
and defined by
\begin{displaymath}
  \V{\lambda}^* =
  \begin{bmatrix}
    676 \\ 196
  \end{bmatrix}
  \qtext{and}
  \M{X}^* = 
  \begin{bmatrix*}[r]
    0 & 3/\sqrt{14} \\
    1 / \sqrt{26} & 2 / \sqrt{14} \\
    -5 / \sqrt{26} & -1/\sqrt{14}
  \end{bmatrix*} \approx
  \begin{bmatrix*}[r]
   0.00 &    0.80 \\
   0.20 &    0.53 \\
  -0.98 &   -0.27 \\
  \end{bmatrix*}.
\end{displaymath}
The matrix $\M{X}^*$ is not orthogonal but is full column rank. We
apply algorithm \Alg{wostd} one hundred times.  We do not apply
randomization (i.e., $\M{V} = \M{I}_n$).  For every run, the predicted
rank is 2 and the solution score is 1 (perfect match). The average
number of attempts (i.e., choosing a random set of $\gamma$-values) to
find a p.s.d.\@ $\M{C}$ is 1.26, and the maximum is 4.

\subsection{Random Non-orthogonal Examples}
In this case, we generate tensors such that $\M{X}^*$ comes from a
matrix with entries drawn from the standard norm distribution whose
columns are normalized (i.e., $\|\MC{x}{k}\|_2=1$ for $k=1,\dots,p$)
and $\M{\lambda}^*$ is the all ones vector.
We apply \Alg{wostd} with no randomization ($\M{V}=\M{I}_n$).
The method fails is it cannot find a p.s.d.\@ $\M{C}$
after 100 attempts.
The results are shown in \Tab{northog}.

In the case of no noise ($\eta=0$), the method is surprising
effective.
Every problem is solved exactly for the even-order tensors ($m=4$ and
$m=6$), which are constructed so that they are positive definite,
guaranteeing that a p.s.d.\@ $\M{C}$ exists.
For the odd-ordered tensor ($m=3$), the method is able to find a p.s.d.\@
$\M{C}$ for 77 out of 100 instances. When it successfully finds the
transformation, the problem is solved exactly.

For the noisy case ($\eta=0.01$), the impact is dramatic. 
In the smallest example ($m=3,n=4,p=2$), only 32 instances can find a
p.s.d.\@ $\M{C}$ matrix, and only 10 instances have a solution score of
0.99 or higher.
For the case  $m=4,n=25,p=3$, the algorithm fails to find a p.s.d.\@
$\M{C}$ in every instance.
For the case $m=6,n=6,p=4$, the algorithm find a p.s.d.\@ $\M{C}$ in a
handful of instances but ultimately fails to solve the problem. We
hypothesize that the problem stems from the fact that the noisy
version of tensor has a rank that is higher than $p$, so the $\M{X}$
that corresponds to the noisy tensor does not have full column rank.

\begin{table}[htbp]
  \centering
    \begin{tabular}{|c@{\,\,\,}c@{\,\,\,}c| %
        >{\color{color1}}r@{~~}  %
        >{\color{color2}}r@{~~}
        >{\color{color3}}r|
        >{\color{color1}}r@{~~}  %
        >{\color{color2}}r@{~~}
        >{\color{color3}}r|
        >{\color{color1}}r@{~~}  %
        >{\color{color2}}r@{~~}
        >{\color{color3}}r|
        >{\color{color1}}r@{~~}  %
        >{\color{color2}}r@{~~}
        >{\color{color3}}r|}
    \hline
    \multicolumn{3}{|c|}{Size} & 
    \multicolumn{12}{c|}{No noise $\eta=0$} \\
    $m$ & $n$ & $p$ &
    \multicolumn{3}{c|}{p.s.d.\@ $\M{C}$?} &
    \multicolumn{3}{c|}{Rank $=p$} &
    \multicolumn{3}{c|}{Rel.~Err.\@ $\leq 10^{-10}$} & 
    \multicolumn{3}{c|}{Soln.~Sc. $\geq 0.99$} \\ 
    \hline
3 &  4 & 2 &  700 &  77 & 11.6 &  700 &  77 & 2.0 &  700 &  77 & 0.0 &  700 &  77 & 1.0 \\
4 & 25 & 3 & 1000 & 100 &  2.5 & 1000 & 100 & 3.0 & 1000 & 100 & 0.0 & 1000 & 100 & 1.0 \\
6 &  6 & 4 & 1000 & 100 &  4.4 & 1000 & 100 & 4.0 & 1000 & 100 & 0.0 & 1000 & 100 & 1.0 \\
    \hline
    \multicolumn{3}{|c|}{Size} & 
    \multicolumn{12}{c|}{Noise $\eta=0.01$} \\
    $m$ & $n$ & $p$ &
    \multicolumn{3}{c|}{p.s.d.\@ $\M{C}$?} &
    \multicolumn{3}{c|}{Rank $=p$} &
    \multicolumn{3}{c|}{Rel.~Err.\@ $\leq$ 0.1} & 
    \multicolumn{3}{c|}{Soln.~Sc. $\geq 0.99$} \\ 
    \hline
3 &  4 & 2 &  238 &  32 & 20.1 &    0 &   0 & 4.0 &   76 &  22 &  0.4 &   17 &  10 & 0.8 \\
4 & 25 & 3 &    0 &   0 &  --- &    0 &   0 & --- &    0 &   0 &  --- &    0 &   0 & --- \\
6 &  6 & 4 &  256 &  35 & 16.9 &    0 &   0 & 6.0 &    0 &   0 & 16.8 &    0 &   0 & 0.4 \\
\hline
    \end{tabular}
  \caption{Random non-orthogonal examples. We create 100 random instances
    and run the method 10 times per instance for a total of 1000
    runs. For each metric, we report \textcolor{color1}{the total number of
      runs where the metric meets the desired criteria}, 
    \textcolor{color2}{the number of instances where at least one run
      meets the desired criteria}, and
    \textcolor{color3}{the mean value of the metric for all
      successful attempts}. 
    For the p.s.d.\@
    $\M{C}$ column, the final number is
      the mean number of attempts needed to find a p.s.d.\@ $\M{C}$
      whenever it occurs.}
  \label{tab:northog}
\end{table}

\section{Conclusions}
If a symmetric tensor $\T{A}$ is known to have a factor matrix $\M{X}$
with orthogonal columns, 
we show that it is possible to solve the symmetric orthogonal tensor
decomposition algorithm via a straightforward matrix eigenproblem.
The method appears to be effective even in the presence of a small
amount of noise.
This is an improvement over previous work \cite{AnGeHsKa12} that
proposed solving the problem iteratively using the tensor eigenvalue
power method and deflation.

We also consider the application of whitening as proposed by
\cite{AnGeHsKa12} in the case where a symmetric tensor $\T{A}$ is
known to have a factor matrix $\M{X}$ that is full rank, but does not
have orthogonal columns. In the noise-free case, the methods works
extremely well. In the case that a small amount of noise is added,
however, the current method is much less effective. Improving
performance in that regime is a potential topic of future study.

\subsection*{Acknowledgments} 
I am grateful to Anima Anandkumar (UC Irvine) for motivating this work with her talk and
 at the
Fields Institute Workshop on Optimization and Matrix Methods in Big
Data.
I am indebted to my Sandia colleagues Grey Ballard and Jackson Mayo
for helpful feedback on this manuscript.

\bibliographystyle{siamdoi}


\end{document}